\documentstyle[12pt,amssymb,amscd]{amsart}  

\newtheorem{thm}{Theorem}[subsection]

\newtheorem{cor}[thm]{Corollary}

\newtheorem{conj}[thm]{Conjecture}

\theoremstyle{remark}
\newtheorem{remark}[thm]{Remark}

\theoremstyle{definition}

\numberwithin{equation}{section}

\newcommand{\g}{{\frak g}}

\newcommand{\ad}{\operatorname{ad}}

\thanks
{Supported
in part by NSF grant DMS-9504522.}

\begin{document}

\title{Formality conjectures for chains}

\author{B.Tsygan}
\address{Department of Mathematics, The Pennsylvania State University,
University Park, PA 16802, USA}
\email{tsygan@@math.psu.edu}

\maketitle

\section{Introduction}
In \cite{K}, Kontsevich proved a classification theorem for deformation quantizations of $C^{\infty}(M)$ where $M$ is a smooth manifold. This theorem asserts that the set of isomorphism classes of deformations of $C^{\infty}(M)$ is in one-to-one correspondence with the set of equivalence classes of formal Poisson structures on $M$. This theorem follows from a more general theorem: the differential graded algebra $\g^{\bullet}_S(M)$ of multi-vector fields on $M$ is equivalent to the differential graded algebra  $\g^{\bullet}_G(C^{\infty}(M))$ of Hochschild cochains of $C^{\infty}(M)$ (we recall exact definitions and statements from \cite{K} and references thereof in section \ref{s:foth}). In other words, the algebra of Hochschild cochains is {\it {formal}} (equivalent to its cohomology).

In this paper, we state a formality conjecture about the Hochschild and cyclic chain complexes of the algebra $C^{\infty}(M)$. It is well known that for any algebra $A$ the Hochschild chain complex $C_{\bullet}(A,A)$ and the negative cyclic complex  $CC^-_{\bullet}(A)$  are modules over the Lie algebra of Hochschild cochains $\g^{\bullet}_G(A)$. Therefore, by virtue of the Kontsevich formality theorem, both the Hochschild (resp. negative cyclic) complex and the graded space of differential forms (resp. de Rham complex) are strong homotopy modules over $\g^{\bullet}_G(C^{\infty}(M))$. We conjecture that those modules are equivalent in an appropriate sense (cf. section \ref{s:foc}), or, in other words, that Connes' quasi-isomorphism from \cite{C}, \cite{L} is, in the right sense, $\g^{\bullet}_G(C^{\infty}(M))$-equivariant. As in \cite{K}, the correct language for stating our conjectures is that of homotopical algebra of Stasheff.

We derive several consequences from the above conjecture in section \ref{hocyde}. First, we compute the Hochschild and cyclic homology of deformed algebras given by the classification theorem of Kontsevich. In particular, we compute the space of traces on such a deformed algebra. Then, in subsection \ref{ss:ahat}, we show how to construct the $\widehat{A}$ class of an arbitrary Poisson manifold. In the case of a regular Poisson structure, this class is, conjecturally, the $\widehat{A}$ class of the tangent bundle to the foliation of symplectic leaves.

Finally, in section \ref{s:ga} we outline a possible proof of our conjectures, as well as their generalization, along the lines of a recent work of Tamarkin \cite{T}.

I am thankful to P. Bressler, G. Halbout, M. Kontsevich, R. Nest, J. Stasheff, D. Tamarkin, and S. Voronov for helpful discussions.
\section{Formality theorem of Kontsevich}  \label{s:foth}
\subsection{Classification of star products}  \label{ss:clasta}

Let $A_0$ be an associative unital algebra over a commutative unital ring $k$. {\it {A deformation}} \cite{G} of $A_0$ is a formal power series 
$$a*b=ab+\sum_{m=1}^{\infty}t^m P_m(a,b)$$
where $P_m: A_0\times A_0\rightarrow A_0$ are $k$-bilinear forms such that

a) The product $*$ is associative

b) $1*f=f*1=f, \; \; f\in A_0$.

{\it {An isomorphism}} of two deformations $*$, $*'$ is a formal power series $T(a)=a+\sum_{m=1}^{\infty}t^m T_m(a)$ such $T(a*b)=T(a)*'T(b)$ for $a,b$ in $A_0$.

Let $M$ be a $C^{\infty}$ manifold. {\it {A deformation quantization}} of $C^{\infty}(M)$, or {\it {a star product}}, is a deformation of $A_0= C^{\infty}(M)$ such that $P_m$ are bidifferential operators \cite{BFFLS}. {\it {An isomorphism}} of two star products is an isomorphism of corresponding deformations such that the operators $T_m$ are differential.

Given a star product on a $C^{\infty}$ manifold $M$, one defines a Poisson bracket on $C^{\infty}(M)$ by
\begin{equation}
\{f,g\}=P_1(f,g)-P_1(g,f)  \label{eq:pbr}
\end{equation}
For a bivector field $\pi$, put 
\begin{equation}
\{f,g\}_{\pi}=<\pi, df \wedge dg>  \label{eq:pbrpi}
\end{equation}
 It follows from associativity of $*$ modulo $t^2$ that the Poisson bracket (\ref{eq:pbr}) is necessarily of the form $\{f,g\}_{\pi_0}$ for some bivector field $\pi_0$. 

Recall that for a bivector field $\pi$ there exists unique trivector field $[\pi,\pi]_S$ such that 
\begin{equation} \label{eq:bras}
\{f,\{g,h\}_{\pi}\}_{\pi} + \{g,\{h,f\}_{\pi}\}_{\pi}+\{h,\{f,g\}_{\pi}\}_{\pi} = <[\pi,\pi]_S, df \wedge dg \wedge dh>
\end{equation}
for any smooth functions $f$, $g$, and $h$. The expression $[\pi,\pi]_S$ is quadratic in $\pi$. The polarization of $[\pi,\pi]_S$ is a symmetric bilinear form $[\pi,\psi]_S$ with values in the space $\Gamma(M,\wedge^3T)$ of trivector fields. A bivector field  $\pi$ such that $[\pi,\pi]_S=0$ is by definition {\it {a Poisson structure}} on $M$.  It follows from associativity of $*$ modulo $t^3$ that for any star product $*$, the bivector field $\pi_0$ (formula (\ref{eq:pbrpi})) is a Poisson structure.

{\it {A formal Poisson structure}} is by definition a formal power series $\pi = \sum_{m=0}^{\infty}t^{m+1} \pi _m$ such that $[\pi,\pi]_S = 0$ in $\Gamma(M,\wedge^3T)[[t]]$.Every Poisson structure $\pi$ defines a formal Poisson structure $t\pi$. Two formal Poisson structures $\pi$ and $\pi '$ are {\it {equivalent}} if there is a formal power series $X=\sum_{m=1}^{\infty}t^m X_m$ such that $\pi ' = \exp (L_X) \pi$.
\begin{thm} \label{thm:konts}
(Kontsevich, \cite{K}). There is a bijection, natural with respect to diffeomorphisms, between the set of equivalence classes of formal Poisson structures on $M$ and the set of isomorphism classes of deformation quantizations of $C^{\infty}(M)$.
\end{thm}

If  $\pi = \sum_{m=0}^{\infty}t^{m+1} \pi _m$ is a formal Poisson structure, we will denote by $*_{\pi}$ a star product from the equivalence class corresponding to $\pi$ by the above theorem. The Poisson structure associated to $*_{\pi}$ by formulas (\ref{eq:pbr}) and (\ref{eq:pbrpi}) will then be equal to $\pi_0$.
\subsection{Hochschild cochains}  \label{ss:coch}
For a unital algebra $A$ over a commutative unital ring $k$, for $n\geq 0$ let 
\begin{equation}  \label{eq:hoco}
\tilde{C}^n(A,A) = Hom (A^{\otimes n}, A)
\end{equation}
If $A=C^{\infty}(M)$, we require that $\tilde{C}^n(A,A)$ consist of those maps 
from $A^{\otimes n}$ to $A$ which are multi-differential.

Define {\it {the Gerstenhaber bracket}} (cf. \cite{G})
$$[\;,\;]_S: \tilde{C}^n(A,A) \otimes \tilde{C}^m(A,A)  \rightarrow \tilde{C}^{n+m-1}(A,A)$$
as follows. For $D \in \tilde{C}^n(A,A)$ and $E \in \tilde{C}^m(A,A)$ put
\begin{equation} \label{eq:circ}
\gathered
	(D\circ E)(a_1,\dots,a_{n+m-1})=\\
=\sum_{j = 0}^{n-1}
	(-1)^{(m-1)j}
	D(a_1,\dots,a_j,  
     E(a_{j+1},\dots,a_{j+m}),\dots);
\endgathered
\end{equation}

\begin{equation} \label{eq:lbra}
        [D, \; E]_S= D\circ E - (-1)^{(n-1)(m-1)}E\circ D
\end{equation}
The bracket $[\;,\;]_S$ turns $\tilde{C}^{\bullet +1}(A,A)$ into a graded Lie algebra (\cite{G}). Put 
\begin{equation} \label{eq:m}
m(a,b)=ab
\end{equation}
for $a,b \in A$.
One has $[m,m]=0$ (this is equivalent to $m$ being associative), so the operator \begin{equation} \label{eq:del}
\delta = [m,?]: \tilde{C}^{\bullet}(A,A) \rightarrow \tilde{C}^{\bullet +1}(A,A)
\end{equation}
satisfies $\delta ^2 = 0$. The complex $(\tilde{C}^{\bullet}(A,A), \delta)$ is called the unnormalized Hochschild cochain complex of $A$ with coefficients in the bimodule $A$. The cohomology of this complex is denoted by $H^\bullet (A,A)$ (the Hochschild cohomology).

Define the (normalized) {\it {Hochschild cochain complex}} of $A$ with coefficients in $A$ by
\begin{equation}  \label{eq:hocon}
{C}^n(A,A) = Hom (\overline{A}^{\otimes n}, A)
\end{equation}
where $\overline{A} = A / k 1$. It is easy to see that $C^{\bullet} (A,A)$ is a subcomplex and a graded Lie subalgebra of $\tilde{C}^{\bullet} (A,A)$. It is well known that the embedding of $C^{\bullet} (A,A)$ into $\tilde{C}^{\bullet }(A,A)$ is a quasi-isomorphism (\cite{CE}).

Put 
\begin{equation} \label{eq:allr}
\g^{\bullet}_G(A) = C^{\bullet + 1}(A,A)
\end{equation}
The differential $\delta$ and the bracket $[,]_S$ make $\g^{\bullet}_G(A)$ a differential graded Lie algebra.

\begin{remark} \label{remark:cup}
Formula
\begin{equation} \label{eq:cup}
\gathered
(D \smile E)(a_1, \dots, a_{n+m}) = \\
=(-1)^{nm}D(a_1, \dots, a_n)E(a_{n+1}, \dots, a_{n+m})
\endgathered
\end{equation} 
defines an associative product on $C^{\bullet}(A,A)$. This product is compatible with the differential $\delta$, therefore $C^{\bullet}(A,A)$ is a differential graded algebra.
\end{remark}
The following theorem is essentially contained in \cite{HKR}
\begin{thm} \label{thm:quisco}
The formula
$$ D_{\pi}(a_1, \dots, a_n) = <\pi, da_1 \dots da_n>$$
defines a quasi-isomorphism 
$$(\Gamma(T, \wedge^{\bullet} T),0) \rightarrow C^{\bullet}(C^{\infty}(M), C^{\infty}(M))$$ Under the isomorphism 
$$\Gamma(T, \wedge^{\bullet} T) \rightarrow H^{\bullet}(C^{\infty}(M), C^{\infty}(M))$$
which is induced by this map on cohomology, the bracket induced by $[\;,\;]_G$ becomes the Schouten bracket $[\;,\;]_S$ (\ref{eq:sch}) and the product induced by the cup product becomes the wedge product.
\end{thm}
\subsection{$L_{\infty}$ algebras}  \label{ss:lin}
{\it {An $L_{\infty}$ algebra}} is a graded vector space $\g ^{\bullet}$ together with a coderivation $\partial$ of the free cocommutative coalgebra $S(\g ^{\bullet}[1])$ such that degree of $\partial$ is $1$ and $\partial ^2 = 0$. 

Put 
\begin{equation} \label{eq:vnesh}
\wedge^{m}(\g^{\bullet})=S^{m}(\g^{\bullet}[1])[-m]
\end{equation}
One has
\begin{equation} \label{eq:vnesh1}
\wedge^{\bullet}(\g^{\bullet})=T(\g^{\bullet})/I
\end{equation}
where $I$ is the two-sided ideal generated by all the elements $$xy-(-1)^{(|x|-1)(|y|-1)}yx$$ where $x, \;y$ are homogeneous elements of $\g^{\bullet}$.
A coderivation $\partial$ of degree $1$ is uniquely determined by its composition with the projection $S(\g^{\bullet}[1]) \rightarrow \g^{\bullet}[1]$ which is a sequence of maps
\begin{equation} \label{eq:holie}
[,\dots,]_n : \wedge^n \g^{\bullet} \rightarrow \g^{\bullet}
\end{equation}
of degree $2-n$, $n=1,2,\dots$
The condition $\partial ^2 = 0$ is equivalent to the following: for any homogeneous elements $x_{1}, \dots, x_{n}$ of $\g^{\bullet}$,
\begin{equation} \label{eq:holie1}
\sum \pm [[x_{i_1}, \dots, x_{i_p}]_p, x_{j_1}, \dots, x_{j_q}]_{q+1} = 0
\end{equation}
where the sum is taken over all ${i_1} < \dots < {i_p}$, ${j_1} < \dots < {j_q}$ such that $\{1, \dots, n\}$ is the disjoint union of $\{{i_1}, \dots, {i_p}\}$ and $\{{j_1} , \dots ,{j_q}\}$. The signs $\pm$ are computed by the following rule: whenever a transposition of $x$ and $y$ appears, the result is multiplied by the sign $(-1)^{(|x|-1)(|y|-1)}$.

In particular, $\delta x = [x]_1$ is a differential on $\g ^{\bullet}$ and the bracket $[\;,\;]_2$ induces a graded Lie algebra structure on the cohomology of the complex $(\g ^{\bullet}, \delta)$.

Any differential graded algebra $(\g ^{\bullet}, [\;,\;], \delta)$ is an $L_{\infty}$ algebra if one puts $[\;]_1 = \delta$, $[\;,\;]_2 = [\;,\;]$ and $[\;,\dots,\;]_n=0$ for $n>2$.

Given two $L_{\infty}$ algebras $\g ^{\bullet}$ and ${\frak{h}} ^{\bullet}$, {\it {an $L_{\infty}$ morphism}} $f: \g ^{\bullet} \rightarrow {\frak{h}} ^{\bullet}$ is a morphism of differential graded coalgebras $S(\g ^{\bullet}[1]) \rightarrow S({\frak{h}} ^{\bullet}[1])$. A morphism of graded coalgebras, without assuming that it commutes with differentials, is uniquely determined by its composition with the projection $S({\frak{h}} ^{\bullet}[1]) \rightarrow {\frak{h} }^{\bullet}[1]$, which is a sequence of linear maps $f_n : \wedge ^n \g ^{\bullet} \rightarrow {\frak{h}} ^{\bullet}$ of degree $1-n$. The condition that these maps define a morphism of differential coalgebras is equivalent to the following: for any homogeneous elements $x_1, \dots, x_n$ of $\g^{\bullet}$,
\begin{equation} \label{eq:molal}
\gathered
\sum \pm f_{q+1}([x_{i_1}, \dots, x_{i_p}]_p, x_{j_1}, \dots, x_{j_q}) = \\
\sum \pm \frac{1}{k!}[f_{n_1}(x_{i_{11}}, \dots, x_{i_{1n_1}}), \dots,f_{n_k}(x_{i_{k1}}, \dots, x_{i_{kn_k}})]_k
\endgathered
\end{equation}
The signs $\pm$ in (\ref{eq:molal}), and in the sum in the left hand side, are as in (\ref{eq:holie}). The sum in the right hand side is taken over all $k \geq 1$ and all ${i_{r1}}< \dots <{i_{rn_r}}$, $1\leq r \leq k$, such that $\{1,\dots, n\}$ is a disjoint union of $\{{i_{r1}}, \dots ,{i_{rn_r}}\}$, $1\leq r \leq k$.

In particular, $f_1$ is a morphism of complexes. We say that f is an $L_{\infty}$ quasi-isomorphism if $f_1$ is a quasi-isomorphism.
\subsection{Formality theorem}  \label{ss:fothm}
Recall that for a manifold $M$ one can define the Schouten-Nijenhujs bracket 
\begin{equation} \label{eq:sch}
[\;,\;]_S: \Gamma(M, \wedge^n T) \otimes \Gamma(M, \wedge^m T) \rightarrow \Gamma(M, \wedge^{n+m-1} T)
\end{equation} as the unique bilinear operation satisfying the following conditions:
\begin{enumerate}
\item for $X \in \Gamma(M, T)$, $[X, \pi]_S = L_X \pi$
\item for $f,g \in \Gamma(M, \wedge^0 T)$, $[f,g]_S = 0$
\item the bracket $[\;,\;]_S$ turns $\Gamma(M, \wedge^{\bullet +1} T)$ into a graded Lie algebra
\item for any $\pi, \varphi, \psi \in \Gamma(M, \wedge^{\bullet} T)$,
$$[\pi, \varphi \wedge \psi]_S = [\pi, \varphi]_S \wedge \psi + (-1)^{|\pi|(|\varphi|+1)}\varphi \wedge [\pi, \psi]_S$$
\end{enumerate}
(for $\pi \in \Gamma(M, \wedge^n T)$, we write $|\pi|=n-1$). When $n=m=2$, the above bracket coincides with the polarized bracket from (\ref{eq:bras}).

Denote by $\g^{\bullet}_S(M)$ the differential graded Lie algebra $\Gamma(M, \wedge^{\bullet +1} T)$ with the bracket $[\;,\;]_S$ and the differential $\delta = 0$.
\begin{thm} \label{thm:fothm}
(Kontsevich, \cite{K}). There exists natural $L_{\infty}$ quasi-isomorphism
$$K: \g_S^{\bullet}(M) \rightarrow \g_G^{\bullet}(C^{\infty}(M))$$
The component $K_1$ of $K$ coincides with the quasi-isomorphism from Theorem \ref{thm:quisco}
\end{thm}

{\bf {Proof of Theorem \ref{thm:konts}}.} For any differential graded Lie algebra $(\g ^{\bullet}, [\;,\;], \delta)$, define 
\begin{equation} \label{eq:mc}
MC(\g ^{\bullet}) = \{ \pi \in t\g ^1[[t]]\;\; |\;\; \delta \pi +\frac{1}{2} [\pi, \pi] = 0\}
\end{equation}
The group $\exp(t\g ^0[[t]])$ acts on $MC(\g^{\bullet})$ by
$$\delta + \exp(X)\pi = \exp(\ad (X))(\delta + \pi)$$
Put 
\begin{equation} \label{eq:mod}
M(\g ^{\bullet})=MC(\g ^{\bullet}) / \exp(t\g ^{0} [[t]])
\end{equation}
For an $L_{\infty}$ homomorphism $f: \g ^{\bullet} \rightarrow
 {\frak{h}} ^{\bullet}$, there is a well-defined  map 
$$M(\g ^{\bullet}) \rightarrow M({\frak{h}} ^{\bullet})$$
 induced by
\begin{equation} \label{eq:mapmod}
\pi \mapsto \sum_{n=1}^{\infty} \frac{1}{n!}f_n(\pi, \dots, \pi)
\end{equation}
If $f$ is an $L_{\infty}$ quasi-isomorphism then the above map is a bijection. Finally, note that $M(\g^{\bullet}_S(M))$ is the set of equivalence classes of formal Poisson structures on $M$ and $M(\g ^{\bullet}_G(A_0))$ is the set of isomorphism classes of deformations of $A_0$ for any algebra $A_0$. Indeed, the equation $\delta \Pi + \frac{1}{2}[\Pi, \Pi] = 0$ is equivalent to $[m+\Pi, m+\Pi]=0$, which is equivalent to the product $a*b = ab + \Pi(a,b)$ being associative.
\begin{remark} \label{rmk:modlinf}
One can easily define spaces $MC(\g ^{\bullet})$ and $M(\g ^{\bullet})$ for any $L_{\infty}$ algebra $\g ^{\bullet}$. For example, the Maurer-Cartan equation from (\ref{eq:mc}) becomes 
$$\sum_{n=1}^{\infty}\frac{1}{n!}[\pi, \dots, \pi]_n = 0$$
\end{remark}
\section{Formality conjectures for Hochschild and cyclic chains} \label{s:foc}
\subsection{Hochschild and cyclic chain complexes}  \label{ss:hocy}
For an algebra $A$ over $k$, define for $n \geq 0$
$$C_n(A,A)=A \otimes \overline{A} ^{\otimes n}$$
(recall that $\overline{A} = A / k1$);
define $b:C_n(A,A) \rightarrow C_{n-1}(A,A)$ by
\begin{equation} \label{eq:b}
\gathered
b(a_0 \otimes \dots \otimes a_n) = (-1)^n a_na_0 \otimes a_1 \otimes \dots \otimes a_{n-1} + \\
\sum _{i=0}^{n-1}(-1)^i a_0 \otimes \dots \otimes a_ia_{i+1} \otimes \dots a_n
\endgathered
\end{equation}
One has $b^2 = 0$. The complex $(C_{\bullet}(A,A), b)$ is called the Hochschild chain complex of $A$ with coefficients in the bimodule $A$. The homology of this complex is denoted by $H_{\bullet}(A,A)$, or $HH_{\bullet}(A)$, and is called {\it {the Hochschild homology}} of $A$.
\begin{remark} \label{rmk:tens}
If $A=C^{\infty}(M)$, one has to use one of the following three definitions of tensor powers of $A$:
\begin{enumerate}
\item $A^{\otimes n+1} = C^{\infty}(M^{n+1})$
\item $A^{\otimes n+1} = \text{germs}_{\Delta}C^{\infty}(M^{n+1})$
\item $A^{\otimes n+1} = \text{jets}_{\Delta}C^{\infty}(M^{n+1})$
\end{enumerate}
where $\Delta$ is the diagonal in $M^{n+1}$.
One defines $A\otimes \overline{A}^{\otimes n}$ accordingly. All the definitions above lead to the same answer for the Hochschild cohomology: the map
\begin{equation} \label{eq:mu}
\mu:( C_{\bullet} (C^{\infty}(M), C^{\infty}(M)),b) \rightarrow (\Omega^{\bullet}(M), 0) 
\end{equation}
defined by 
\begin{equation} \label{eq:mu1}
\mu (a_0 \otimes \dots \otimes a_n)=\frac{1}{n!}a_0 da_1 \dots da_n
\end{equation}
is a quasi-isomorphism of complexes (cf. \cite{C}).
\end{remark}

For $D \in C^d(A,A)$ and $a_0, \dots, a_n \in A$, define
\begin{equation} \label{eq:ld}
\gathered
L_D(a_0, \dots, a_n)=\sum_{i=0}^{n-d}(-1)^{(d-1)(i+1)}a_0 \otimes \dots a_i \otimes D(a_{i+1},\dots, a_{i+d}) \otimes \dots a_n + \\
\sum_{j=n-d}^{n}(-1)^{n(j+1)}D(a_{j+1}, \dots, a_0, \dots) \otimes a_{d+j-n} \dots \otimes a_j
\endgathered
\end{equation}
One has 
$$ [L_D, L_E]=L_{[D,E]_G}$$
and
$$b=L_m$$
Thus the operators $L_D$ define an action of the differential graded Lie algebra $\g _G ^{\bullet}(A)$ on the complex $C_{\bullet}(A,A)$.

Recall that {\it {the cyclic differential}} 
\begin{equation}  \label{eq:bcap}
B: C_n(A,A) \rightarrow C_{n+1}(A,A)
\end{equation}
is defined by 
\begin{equation}  \label{eq:bcap1}
B(a_0 \otimes \dots a_n) = \sum_{i=0}^{n}(-1)^{ni}1\otimes a_i \otimes \dots \otimes a_0 \otimes \dots \otimes a_{i-1}
\end{equation}
One has $b^2 = Bb+bB = B^2 = 0$, as well as 
$$[B, L_D]=0$$
Following Getzler's notation, put 
\begin{equation} \label{eq:ccmin}
CC^-_{\bullet}(A) = (C_{\bullet}(A,A)[[u]],b+uB)
\end{equation}
where $u$ is a formal variable of degree $-2$. One sees that $CC^-_{\bullet}(A)$
is a differential graded module over the differential graded algebra $\g ^{\bullet}_G (A)$ for any algebra $A$.

Now consider the space $\Omega^{\bullet}(M)$ of differential forms on a manifold $M$. For a multivector field $\pi \in \Gamma (M, \wedge ^d T)$,
put
\begin{equation} \label{eq:lpi}
L_{\pi}=[d,i_{\pi}]
\end{equation}
 where $i_{\pi}$ is the contraction 
\begin{equation} \label{eq:ipi}
i_{\pi}: \Omega^{\bullet}(M) \rightarrow \Omega^{\bullet - d}(M)
\end{equation}
It is well known that 
\begin{equation} \label{eq:lipi}
L_{[\pi,\psi]_S}=[L_{\pi},L_{\psi}]
\end{equation}
Therefore operators $L_{\pi}$ define an action of the algebra $\g_S^{\bullet}(M)$ on the graded space $\Omega^{\bullet}(M)$, as well as on the complex $(\Omega^{\bullet}(M)[[u]], ud)$.
\begin{thm} \label{thm:mucyc}
(Connes, \cite{C}). The formula
$$\mu(a_0 \otimes \dots \otimes a_n) = \frac{1}{n!}a_0 da_1 \dots d a_n$$
defines a functorial ${\Bbb{C}}[[u]]$-linear quasi-isomorphism
\begin{equation} \label{eq:mucyc}
CC_{\bullet}^-(C^{\infty}(M)) \rightarrow (\Omega^{\bullet}(M)[[u]],ud)
\end{equation}
\end{thm}
Under the homomorphisms
\begin{equation} \label{eq:klipi}
H_{\bullet}(C^{\infty}(M),C^{\infty}(M)) \rightarrow \Omega^{\bullet}(M)
\end{equation}
\begin{equation} \label{eq:piipi}
HC_{\bullet}^-(C^{\infty}(M)) \rightarrow H^{\bullet}(M)[[u]]
\end{equation}
which are induced by $\mu$ on cohomology, the operators induced by $L_D$ become $L_{\pi}$ where $D$ is a Hochschild cocycle and $\pi$ is the image of the cohomology class of $D$ under the isomorphism from Theorem \ref{thm:quisco}.

\subsection{$L_{\infty}$ modules}  \label{ss:limo}
 {\it {An $L_{\infty}$ module}} over an $L_{\infty}$ algebra $\g^{\bullet}$ is a graded vector space $M^{\bullet}$ together with a coderivation $\partial_M$ of degree $1$ of the free differential graded comodule $S(\g^{\bullet} [1])\otimes M^{\bullet}$ such that $\partial_M ^2 = 0$. A coderivation $\partial_M$, without the condition $\partial_M ^2 = 0$, is uniquely determined by its composition with the projection of  $S(\g^{\bullet} [1])\otimes M^{\bullet}$ onto $M^{\bullet}$. This composition is a sequence of linear maps 
\begin{equation} \label{eq:fi}
\phi _n : \wedge^n \g^{\bullet} \otimes M^{\bullet} \rightarrow M^{\bullet}
\end{equation}
of degree $1-n$. The condition $\partial_M ^2 =0$ is equivalent to the following: for any homogeneous elements $x_1, \dots, x_n$ of $\g^{\bullet}$ and $m$ of $M^{\bullet}$, 
\begin{equation} \label{eq:linfmodu}
\gathered
\sum \pm \phi_{p+1}(x_{i_1}, \dots, x_{i_p},\phi_q( x_{j_1}, \dots, x_{j_q},m))+\\
+\sum \pm \phi_{q+1}([x_{i_1}, \dots, x_{i_p}]_p, x_{j_1}, \dots, x_{j_q},m) = 0
\endgathered
\end{equation}
where the signs are computed and the sum is taken as in (\ref{eq:holie}).
In particular, 
$$\delta_M = \phi_0$$ is a differential on $M^{\bullet}$.
The maps $\phi _n$, $n \geq 0$, define an $L_{\infty}$ module structure on $M^{\bullet}$ if and only if the maps $\phi _n$, $n\geq 1$, define an $L_{\infty}$ morphism $\g^{\bullet} \rightarrow \text{Hom}(M^{\bullet},M^{\bullet})$ where the right hand side is viewed as a differential graded Lie algebra with the differential $[\delta_M, ?]$.

A morphism of $L_{\infty}$ modules over $\g^{\bullet}, \;$ $\varphi : M^{\bullet} \rightarrow N^{\bullet}$ is by definition a morphism of differential graded comodules  $S(\g^{\bullet} [1])\otimes M^{\bullet} \rightarrow  S(\g^{\bullet} [1])\otimes N^{\bullet}$. Such a morphism is uniquely determined by maps
$$\varphi _n : \wedge ^n \g^{\bullet} \otimes M^{\bullet} \rightarrow N^{\bullet}$$ of degree $-n$, $n \geq 0$, satisfying
\begin{equation} \label{eq:linfimohomo}
\gathered
\sum \pm \varphi _{q+1}([x_{i_1}, \dots, x_{i_p}]_p, x_{j_1}, \dots, x_{j_q},m) +\\
+\sum \pm \varphi_{p+1}(x_{i_1}, \dots, x_{i_p},\phi_q( x_{j_1}, \dots, x_{j_q},m))= \\ =
\sum \pm \phi_{p+1}(x_{i_1}, \dots, x_{i_p},\varphi_q( x_{j_1}, \dots, x_{j_q},m))
\endgathered
\end{equation}
\begin{remark} \label{remark:3edef}   
An equivalent definition of an $L_{\infty}$ module structure on $M^{\bullet}$ is the following. Let $\g^{\bullet}$ be an $L_{\infty}$ algebra and $M^{\bullet}$ a graded vector space. On the graded space $\g^{\bullet} \oplus M^{\bullet}$ consider another grading in which $\g^{\bullet}$ is of degree zero and $M^{\bullet}$ is of degree one. Consider an $L_{\infty}$ algebra structure on $\g^{\bullet} \oplus M^{\bullet}$  such that:
\begin{enumerate}
\item $\g^{\bullet}$ is an $L_{\infty}$ subalgebra of $\g^{\bullet} \oplus M^{\bullet}$ 
\item all the operations $[\;,\dots,\;]_n$ are of degree zero with respect to the second grading on $\g^{\bullet} \oplus M^{\bullet}$ 
\item $[m_1, m_2, \dots]_n = 0$ for any $m_1, \; m_2 \in M^{\bullet}$.
\end{enumerate}
Similarly, one can define a morphism of $L_{\infty}$ modules as an $L_{\infty}$ morphism $f: \g^{\bullet} \oplus M^{\bullet} \rightarrow \g^{\bullet} \oplus N^{\bullet}$  which is of degree zero with respect to the second grading and for which $f_n(m_1,  m_2, \dots ) = 0$ for any $m_1, \; m_2 \in M^{\bullet}$.
\end{remark}
\subsection{Formality conjecture for chains}  \label{focha}
 Let $K$ be the $L_{\infty}$ quasi-isomorphism from Theorem \ref{thm:fothm}. Via $K$, the differential graded Lie algebra $\g^{\bullet}_S(M)$ acts on $C_{\bullet}(C^{\infty}(M), C^{\infty}(M))$ and $CC_{\bullet}^-(C^{\infty}(M))$ as on $L_{\infty}$ modules.
\begin{conj} \label{conj:fococh}
There exists a natural quasi-isomorphism of $L_{\infty}$ modules 
$$\phi: C_{\bullet}(C^{\infty}(M),C^{\infty}(M)) \rightarrow (\Omega^{\bullet}(M),0)$$
such that $\phi _0 $ is the quasi-isomorphism $\mu$ of Connes.
\end{conj}
This conjecture extends to the following 
\begin{conj} \label{conj:focycoch}
There exists a natural ${\Bbb{C}}[[u]]$-linear quasi-isomorphism of $L_{\infty}$ modules
$$\phi: CC_{\bullet}^-(C^{\infty}(M)) \rightarrow (\Omega^{\bullet}(M)[[u]],ud)$$ such that $\phi_0$ is the Connes quasi-isomorphism $\mu$.
\end{conj}
\section{Hochschild and cyclic complexes of deformed algebras} \label{hocyde}
Let $\pi$ be a 
formal Poisson structure on a manifold $M$. The isomorphism from Theorem \ref{thm:konts} provides a star product $*_{\pi}$ on $C^{\infty}(M)$. Put 
\begin{equation} \label{eq:baba}
A(\pi) = (C^{\infty}(M), *_{\pi})
\end{equation}
This is an algebra over $k = {\Bbb{C}}[[t]]$. By $A(\pi)^{\otimes _k (n+1)}$ we will denote $C^{\infty}(M)^{\otimes (n+1)}[[t]]$ (cf. Remark \ref{rmk:tens}); similarly for $A(\pi) \otimes \overline{A(\pi)}^{\otimes _k n}$). If Conjecture \ref{conj:focycoch} is true, then one gets
\begin{cor} \label{cor:ccdef}
There exists a quasi-isomorphism
\begin{equation} \label{eq:cdef}
 \mu^{\pi}: C_{\bullet}(A(\pi),A(\pi))\rightarrow (\Omega^{\bullet}(M)[[t]], L_{\pi})
\end{equation}
which extends to a ${\Bbb{C}}[[u]]$-linear quasi-isomorphism
\begin{equation} \label{eq:ccdef}
 \mu^{\pi} :CC^-_{\bullet}(A(\pi)) \rightarrow (\Omega^{\bullet}(M)[[u]][[t]], L_{\pi}+ud)
\end{equation}
\end{cor}
\begin{pf} Let $K$ be the quasi-isomorphism from Theorem \ref{thm:fothm}. One checks that the formula
\begin{equation} \label{eq:kapi}
\mu^{\pi}(c) = \sum _{n=0}^{\infty}\frac{1}{n!}K_{n}(\pi, \dots, \pi , c)
\end{equation} 
defines a quasi-isomorphism of complexes.
\end{pf}
\begin{remark}
The complexes in the right hand side of formulas (\ref{eq:cdef}, \ref{eq:ccdef}) were studied by Brylinski as quasi-classical approximations to the left hand sides. If $\pi = t \pi _0$ where $\pi _ 0$ is a Poisson structure, then filtration by powers of $t$ defines a spectral sequence with the $E^1$ term equal to the right hand side; this spectral sequence converges to the left hand side (\cite{Bryl}). The above Corollary implies that this spectral sequence degenerates at $E^1$
\end{remark}
\begin{cor} \label{cor:tr}
The space of continuous ${\Bbb {C}}[[t]]$-valued traces on $A(\pi)$ is isomorphic to the space of continuous ${\Bbb {C}}[[t]]$-linear ${\Bbb {C}}[[t]]$-valued functionals on $C^{\infty}(M)$ which annihilate all Poisson brackets $\{f,g\}_{\pi}$\end{cor}
(compare with \cite{CFS}, \cite{Fe}, \cite{NT3} for the symplectic case).
\subsection{The $\widehat{A}$ class of a Poisson manifold} \label{ss:ahat}
 Consider a $C^{\infty}$ manifold M with a Poisson structure $\pi _0$. Assuming that Conjecture 
\ref{conj:focycoch} is true, define the cohomology class $\widehat{A}(\pi_0)$ in $H^{\text{ev}}(M, {\Bbb{C}}[[t]])$ as folows.

Let $\pi = t \pi _0$. Recall that for any $k$-algebra $A$ the periodic cyclic complex of $A$ is defined by
\begin{equation} \label{eq:ccper}
CC^{per}_{\bullet}(A) = (C_{\bullet}(A,A)[u^{-1}, u]], b+uB)
\end{equation}
Localizing the $L_{\infty}$ quasi-isomorphism $\mu^{\pi}$ from Corollary \ref{cor:ccdef} with respect to $u$, one gets an $L_{\infty}$ quasi-isomorphism 
\begin{equation} \label{eq:ccpdef}
 \mu^{\pi} :CC^{per}_{\bullet}(A(\pi)) \rightarrow (\Omega^{\bullet}(M)[u^{-1},u]][[t]], L_{\pi}+ud)
\end{equation}
Now note that the complex in the right hand side of (\ref{eq:ccpdef}) is isomorphic to the complex $(\Omega^{\bullet}(M)[u^{-1},u]][[t]],ud)$. Indeed, $L_{\pi} = [d, i_{\pi}]$ and the desired isomorphism is given by $\exp(-u^{-1}ti_{\pi})$. One gets a quasi-isomorphism 
\begin{equation} \label{eq:ccpdef0}
 CC^{per}_{\bullet}(A(\pi)) \rightarrow (\Omega^{\bullet}(M)[u^{-1},u]][[t]], ud)
\end{equation}
If one views $1$ as an element of $C_0(A(\pi), A(\pi))$, and thus of $CC^{per}_0(A(\pi))$, then the value of the quasi-isomorphism (\ref{eq:ccpdef0}) at $1$ is an element of degree zero in $H^{\bullet}(M, {\Bbb{C}}[[t]])[u^{-1},u]]$, so it can be viewed as an element $\widehat{A}(\pi _0)$ of $H^{\text{ev}}(M, {\Bbb{C}}[[t]])$. Conjecturally, this class does not depend on $t$.

Consider the situation when $\pi _0$ is a regular Poisson structure. In this case, the symplectic leaves of $\pi_0$ form a foliation ${\cal{F}}$. The tangent bundle $T{\cal{F}}$ of this foliation is an $Sp(2n)$-bundle, and one can reduce its structure group to the maximal compact subgroup $U(n)$. Let $\widehat{A}(T{\cal{F}})$ be the $\widehat{A}$ class of the resulting $U(n)$-bundle.

More generally, suppose that $\pi _0$ comes from a symplectic Lie algebroid $({\cal{E}}, \omega) $ (\cite{MK}, \cite{BB}; cf.below for the definitions). This generality was suggested to us by Weinstein. 
\begin{conj} \label{conj:ahat}
If  $\pi _0$ comes from a symplectic Lie algebroid $({\cal{E}}, \omega)$ then 
$$\widehat{A}(\pi_0) = \widehat{A} ({\cal{E}})$$
\end{conj}
Recall that {\it {a Lie algebroid}} is a vector bundle ${\cal{E}}$ whose sections form a sheaf of Lie algebras, together with a morphism of sheaves of Lie algebras ({\it {the anchor map}} )
\begin{equation} \label{eq:anch}
\rho : \Gamma({\cal{E}}) \rightarrow \Gamma(T)
\end{equation}
such that
\begin{equation} \label{eq:ancho}
 [\xi, f\eta] = L_{\rho(\xi)}(f)\eta + f[\xi, \eta]
\end{equation}
where $\xi, \; \eta$ are local sections of ${\cal{E}}$ and $f$ is a local function.  Recall that for a Lie algebroid ${\cal{E}}$, the de Rham complex is defined:
\begin{equation} \label{eq:dre}
\gathered
^{\cal{E}}\Omega^{\bullet}(M) = \Gamma (M, \wedge ^{\bullet}{\cal{E}}^*)\\
^ {\cal{E}}\Omega^{\bullet}(M) \stackrel{d}{\rightarrow}{^ {\cal{E}}\Omega^{\bullet + 1}(M)} \\
(d\omega)(\xi _1, \dots, \xi _{m+1}) = \sum_{i=1}^{m+1}(-1)^{i-1}\rho(\xi _i)\omega(\xi _1, \dots, \widehat{\xi _i}, \dots, \xi_{m+1}) + \\
\\ +\sum _{i<j}(-1)^{i+j} \omega([\xi_i, \xi_j], \dots, \widehat{\xi _i}, \dots, \widehat{\xi _j}, \dots)
\endgathered
\end{equation}

The algebra of {\it {${\cal{E}}$-differential operators}} $^{\cal{E}} D_M$ is the quotient of the enveloping algebra $U(\Gamma ({\cal{E}})$ by the ideal generated by the elements 
$$  \xi f\eta - L_{\rho(\xi)}(f)\eta - f\xi \eta $$
where $\xi, \; \eta$ are local sections of ${\cal{E}}$ and $f$ is a local function.

{\it {A symplectic Lie algebroid}} is a Lie algebroid ${\cal{E}}$ together with a non-degenerate closed ${\cal{E}}$-form $\omega \in { ^ {\cal{E}}\Omega^{\bullet}(M)}$. We denote by $\pi _0 \in \Gamma (M, \wedge ^2{\cal{E}})$ the image of $\omega$ under the isomorphism $\wedge ^2{\cal{E}}^* \rightarrow \wedge ^2{\cal{E}}$ induced by $\omega$. The bivector field $(\wedge ^2 \rho)(\pi _0)$ is a Poisson structure. We will denote this Poisson structure also by $\pi _0$. 

Let us outline the reasoning for which the above conjecture should be true. Corollary \ref{cor:ccdef} is true for $\pi = t \pi _0$ where $\pi _0$ is a regular Poisson structure or, more generally, when $\pi _0$ is given by a symplectic Lie algebroid $({\cal{E}}, \omega) $.The cohomology of the above complex will be denoted by $^{\cal{E}}H^{\bullet}(M)$.
When ${\cal{F}}$ is a foliation with a leafwise symplectic form $\omega$ and ${\cal{E}}$ is the tangent bundle of this foliation, then ${^ {\cal{E}}\Omega^{\bullet}(M)}$ is the de Rham complex of leafwise forms. The anchor map extends to a morphism of complexes $\Omega^{\bullet}(M) \rightarrow {^ {\cal{E}}\Omega^{\bullet}(M)}$.

In \cite{NT1}, we studied star products on $C^{\infty}(M)$ for which the corresponding Poisson structure $\pi _0$ is given by a symplectic Lie algebroid $({\cal{E}}, \omega)$ and the operators $P_m$ are ${\cal{E}}$-bidifferential (call them {\it {$ {\cal{E}}$-deformations}}). We regard two ${\cal{E}}$-deformations as equivalent if there is an equivalence of star products $T=\text{Id}+\sum_{m=1}^{\infty}t^mT_m$ for which all $T_m$ are ${\cal{E}}$-differential operators. We showed that Fedosov's methods from \cite{F} are applicable in this situation. In particular, to any ${\cal{E}}$-deformation one can associate a characteristic class
\begin{equation}  \label{eq:teta}
\theta \in \frac{1}{t}\omega + {^ {\cal{E}}H^2(M, {\Bbb{C}})}[[t]]
\end{equation}
which defines a bijection between the set of equivalence classes of ${\cal{E}}$-deformations and $\frac{1}{t}\omega + {^{\cal{E}}H^2(M, {\Bbb{C}})}[[t]]$. 

Suppose given an ${\cal{E}}$-deformation $A = (C^{\infty}(M), *)$ with the characteristic class $\theta$. In \cite{NT1} and \cite{BNT}, we constructed a ${\Bbb{C}}[[u]]$-linear trace density map 
\begin{equation} \label{eq:muh}
\mu^t: CC_{\bullet}^-(A) \rightarrow (^ {\cal{E}}\Omega^{2n-\bullet}(M)((t))[[u]],d)
\end{equation}
whose localization with respect to $u$ provides a ${\Bbb{C}}[u^{-1},u]]$-linear morphism
\begin{equation} \label{eq:muhp}
\mu^t: CC_{\bullet}^{per}(A) \rightarrow (^ {\cal{E}}\Omega^{2n-\bullet}(M)((t))[u^{-1},u]],d)
\end{equation}
We proved
\begin{thm} \label{thm:rrbnt}
$$\mu ^t (1) = \sum_{p\geq 0}\widehat{A}({\cal{E}})_{2p}u^{n-p}$$
\end{thm}
Now assume for simplicity that $\theta = \frac{1}{t}\omega$.

Let $N=(-1)^{n}t^{k-n}$ on $^{\cal{E}}\Omega^{k}(M)$. Let $*: ^ {\cal{E}}\Omega^{\bullet}(M) \rightarrow ^ {\cal{E}}\Omega^{2n-\bullet}(M)$ be the symplectic star operator. Consider the sequence of morphisms of complexes \begin{equation} \label{eq:mnogo}
\gathered
(^ {\cal{E}}\Omega^{2n-\bullet}(M)((t))((u)),d) \stackrel{N}{\longrightarrow}( ^ {\cal{E}}\Omega^{2n-\bullet}(M)((t))((u)), td) \\
\stackrel {\exp(ut^{-1}i_{\pi_0})}{\longrightarrow} (^ {\cal{E}}\Omega^{2n-\bullet}(M)((t))((u)), td + uL_{\pi_0}) \\
\stackrel{*}{\longrightarrow} (^ {\cal{E}}\Omega^{\bullet}(M)((t))((u)), tL_{\pi_0}+ud) \\
\stackrel{\exp(-tu^{-1}i_{\pi_0})}{\longrightarrow} (^{\cal{E}}\Omega^{\bullet}(M)((t))((u)), ud)
\endgathered
\end{equation} 
Denote the composition of the above maps by $\nu$.
\begin{conj} \label{conj:whocares}
The quasi-isomorphism (\ref{eq:ccpdef0}), composed with 
\begin{equation}  \label{eq:omeom}
\Omega^{\bullet}(M)((t))((u)) \rightarrow ^ {\cal{E}}\Omega^{\bullet}(M)((t))((u)),
\end{equation}
 is equal to $\nu \mu ^t$.
\end{conj}
Compose $\nu$ with the isomorphism
$$^ {\cal{E}}\Omega^{\bullet}(M)((t))((u)) \rightarrow ^ {\cal{E}}\Omega^{\bullet}(M)((t))((u))$$
which is equal to $u^{n-k}$ on $^{\cal{E}} \Omega ^k$.
\begin{remark}
Note the symmetry between the formal variables $u$ and $t$ in the above calculations.
\end{remark}
 Denote the resulting composition by $\nu _0$. We claim that $\nu _0$ is equal to the operator of multiplication by $\exp(-\frac{\omega}{ut})$. Indeed, one checks that $\nu _0 (1)$ is equal to $\exp(-\frac{\omega}{ut})$ (we use the fact that for any $z$ 
$$\exp(zi_{\pi})\frac{1}{n!}\omega^{n}=z^n\exp(z^{-1}\omega)).$$
But $\nu _0$ is a $C^{\infty}(M)$-linear endomorphism of $^ {\cal{E}}\Omega^{2n-\bullet}(M)((t))((u))$, thus it is the operator of multiplication by $\nu _0 (1)$. From Theorem \ref{thm:rrbnt} we see that 
\begin{equation} \label{eq:numu1}
(\nu_0\mu^t)(1)=u^{n}\sum \widehat{A}({\cal{E}})_{2p}u^{-p}
\end{equation}
Therefore
\begin{equation} \label{eq:numu2}
(\nu_0\mu^t)(1)=\sum \widehat{A}({\cal{E}})_{2p}u^{p}
\end{equation}
Combining the above formula with Conjecture \ref{conj:whocares}, one sees that the composition of maps (\ref{eq:omeom}) and (\ref{eq:ccpdef0}) evaluated at $1$ is equal to $\sum \widehat{A}({\cal{E}})_{2p}u^{p}$.
\section{Homotopy Gerstenhaber algebras and modules}  \label{s:ga}
In this section, we will outline a possible proof of the conjectures above, as well as their generalizations. This proof will follow the lines of the recent proof of the Kontsevich formality theorem, due to Tamarkin \cite{T}.
\subsection{Definitions}  \label{ss:gade}
Recall that a graded space $V^{\bullet}$ is {\it {a Gerstenhaber algebra}} if it is a graded commutative associative algebra, $V^{\bullet}[1]$ is a graded Lie algebra, and the two operations on $V^{\bullet}$ satisfy the Leibnitz identity
\begin{equation}  \label{eq:le}
[a,bc]=[a,b]c+(-1)^{(|a|-1)|b|}b[a,c]
\end{equation}
The Hochschild cohomology $H^{\bullet}(A,A)$ of any associative algebra $A$ is a Gerstenhaber algebra on which the product and the bracket are induced by the cup product and the Gerstenhaber bracket respectively ((\ref{eq:cup}), (\ref{eq:lbra})). 

Let us recall a definition of a $G_{\infty}$ algebra, or a strong homotopy Gerstenhaber algebra. For a graded vector space $V^{\bullet}$, consider the free Lie coalgebra $\text{coLie}(V^{\bullet}[1])$ and the free cocommutative coalgebra $S(\text{coLie}(V^{\bullet}[1]))$. The latter graded space has a structure of a Lie coalgebra which is dual to the Berezin-Kirillov-Kostant Lie algebra structure (in the dual language, $S(\text{Lie}(V^{\bullet}[1]))$ is a Poisson algebra).
%
%
%

A structure of a $G_{\infty}$ algebra on $V^{\bullet}$ is by definition a map $\partial : S(\text{coLie}(V^{\bullet}[1])) \rightarrow S(\text{coLie}(V^{\bullet}[1]))$ of degree $1$ which is a coderivation with respect to both coalgebra structures and such that $\partial ^2 = 0$. Any Gerstenhaber algebra is a $G_{\infty}$ algebra. For any $G_{\infty}$ algebra  $V^{\bullet}$, $V^{\bullet}[1]$ is an $L_{\infty}$ algebra.

For any $G_{\infty}$ algebra $V^{\bullet}$, one can define the cochain complex of coderivations of  $S(\text{coLie}(V^{\bullet}[1]))$, with the differential $[\partial, ?]$. The cohomology of this complex is denoted by $H^{\bullet}(V^{\bullet}, V^{\bullet})$. This is a $G_{\infty}$ analogue of Hochschild cohomology.

Define a $G_{\infty}$ morphism $f:V^{\bullet} \rightarrow W^{\bullet}$ to be a morphism of differential graded Poisson coalgebras  $S(\text{coLie}(V^{\bullet}[1])) \rightarrow S(\text{coLie}(W^{\bullet}[1]))$. We say that a $G_{\infty}$ algebra $V^{\bullet}$ is formal if there is a $G_{\infty}$ quasi-isomorphism 
\begin{equation} \label{eq:forgest}
H^{\bullet}(V^{\bullet}, V^{\bullet}) \rightarrow V^{\bullet}
\end{equation}
A standard argument from homological algebra shows that obstructions to formality of a $G_{\infty}$ algebra $V^{\bullet}$ lie in $H^{\bullet}(V^{\bullet}, V^{\bullet})$.
\subsection{Tamarkin's proof}   \label{ss:tap}
In \cite{T}, Tamarkin proves the following
\begin{thm} \label{thm:tam}
For any associative algebra $A$, the Hochschild cochain complex $C^{\bullet}(A,A)$ has a structure of a $G_{\infty}$ algebra whose underlying $L_{\infty}$ algebra is $\g^{\bullet}_G(A)$
\end{thm}
\begin{thm} \label{thm:tam1}
Let $A={\Bbb{C}}[[x_1, \dots, x_n]]$ or $A=C^{\infty}({\Bbb {R}}^n)$. The obstructions to formality of the $G_{\infty}$ algebra  $C^{\bullet}(A,A)$ are equal to zero.
\end{thm}
This shows that the above algebras are formal as $G_{\infty}$ algebras. From this, using an argument with Gelfand-Fuks cohomology as in \cite{K}, one deduces
\begin{thm}
Let $A=C^{\infty}(M)$. Then $C^{\bullet}(A,A)$ is formal as a $G_{\infty}$ algebra.
\end{thm}
\subsection{Generalized formality conjecture for chains}
Conjecture \ref{conj:fococh} can be generalized along the lines of the previous subsection as follows. First, one can define $G_{\infty}$ modules and their homomorphisms as one did in the $L_{\infty}$ case (following any of the above definitions, for example the one from Remark \ref{remark:3edef}).

The problem with this definition is that, for example, $\Omega ^{\bullet}(M)$ is not a Gerstenhaber module over the Gerstenhaber algebra $\Gamma (\wedge ^{\bullet}(T))$. To correct this, for any Gerstenhaber algebra $V^{\bullet}$ define a new Gerstenhaber algebra $V^{\bullet}[\epsilon]$ by
\begin{equation} \label{epsilon}
(a+\epsilon b)(c+\epsilon d)=ac+\epsilon(bc + (-1)^{|a|}ad + (-1)^{|a|}[a,c])
\end{equation}
\begin{equation} \label{epsilonlie}
[a+\epsilon b,\;c+\epsilon d]=[a,c]+\epsilon([b,c] + (-1)^{|a|+1}[a,d ])
\end{equation}
(This is a deformation of $V^{\bullet}$ with an odd parameter along the Poisson bracket; the specifics of the graded case is that the deformed algebra remains commutative. Note that an isomorphism of this deformation to the trivial one is precisely a BV operator). Using Tamarkin's methods, one can prove that for any algebra $A$ there is a $G_{\infty}$ structure on $C_{\bullet}(A,A)[\epsilon]$ which induces the structure (\ref{epsilon}), (\ref{epsilonlie}) on $H_{\bullet}(A,A)[\epsilon]$.
\begin{conj} \label{conj:ger}
For any associative algebra $A$, the Hochschild chain complex $C_{\bullet}(A,A)$ is a $G_{\infty}$ module over the $G_{\infty}$ algebra 
$C^{\bullet}(A,A)[\epsilon]$. 
The underlying structure of an  $L_{\infty}$ module over  $C^{\bullet}(A,A)$ on $C_{\bullet}(A,A)$ is given by the action of $C^{\bullet}(A,A)$ by operators $L_D$ (\ref{eq:ld}).
\end{conj}
If the above conjecture is true then, by virtue of Theorem \ref{thm:tam}, both 
$C_{\bullet}(C^{\infty}(M),C^{\infty}(M))$ and $\Omega ^{\bullet}(M)$ are $G_{\infty}$ modules over 
$\Gamma(M, \wedge ^{\bullet}T)[\epsilon]$.
\begin{conj} \label{conj:fogerts}
There is a quasi-isomorphism of $G_{\infty}$ modules 
$$C_{\bullet}(C^{\infty}(M),C^{\infty}(M)) \rightarrow  \Omega ^{\bullet}(M)$$
\end{conj}
To generalize Conjecture \ref{conj:focycoch}, first note that the operator $\frac{\partial}{\partial \epsilon}$ is a BV operator on the algebra $V^{\bullet}[\epsilon]$ (\ref{epsilon}), (\ref{epsilonlie}). Conjecturally, in an appropriate sense, $C_{\bullet}^-(A, A)[\epsilon]$ is a homotopy BV algebra and $CC_{\bullet}^-(A)[\epsilon]$ is a homotopy BV module over it.

Let us finish with a partial case of the statement before Conjecture \ref{conj:ger} which can be obtained by explicit computation.
\begin{thm}  \label{thm:tsydal}
\cite{DT}. On the Hochschild chain complex $C_{\bullet}(A,A)$, there is a structure of an $L_{\infty}$ module over  
$\g_G^{\bullet}(A)[\epsilon]$ whose restriction to  $\g_G^{\bullet}(A)$ is given by the operators $L_D$.
\end{thm}

\end{document}